\documentclass[12pt]{article}
\usepackage[cp1251]{inputenc}
\usepackage{amsmath}
\usepackage{amssymb}
\usepackage[english]{babel}
\usepackage[T2A]{fontenc}
\usepackage{epsfig}

\begin{document}

\begin{large}

\centerline{\bf On the flexibility of Siamese dipyramids }

\end{large}

\medskip

\centerline{ I. Fesenko (V.Karazin Kharkiv National University, Kharkiv, Ukraine)}

\centerline{ V. Gorkavyy (B.Verkin Institute for Low Temperature
Physics and Engineering, Kharkiv, Ukraine)}

\date{empty}

\bigskip

\begin{small}

{\it Abstract}. Polyhedra called Siamese dipyramids are known to
be non-flexible, however their physical models behave like
physical models of flexible polyhedra. We discuss a simple
mathematical method for explaining the model flexibility of the
Siamese dipyramids.



\end{small}

\bigskip

\bigskip

{\bf \large Introduction}

\medskip

This paper deals with particular polyhedral surfaces called {\it
Siamese dipyramids}, and it is aimed to explain their model
flexibility.

Siamese dipyramids were described by M. Goldberg in \cite{Goldberg}\footnote{ Actually, in
\cite{Goldberg} the particular case of equilateral pentagonal Siamese dipyramids was discussed
only.}, see also \cite[p.222-225]{Cromwell}. By definition, a general $n$-gonal Siamese dipyramid
consists of $4n$ isosceles triangular faces and looks like two $n$-gonal dipyramids stuck together,
see Fig.1. For any $n\geq 3$ there exists a wide variety of different $n$-gonal Siamese dipyramids
which differ in the lengths of edges as well as in the spatial shapes.

\begin{figure}[ht!]
    \centering
    \epsfig{file=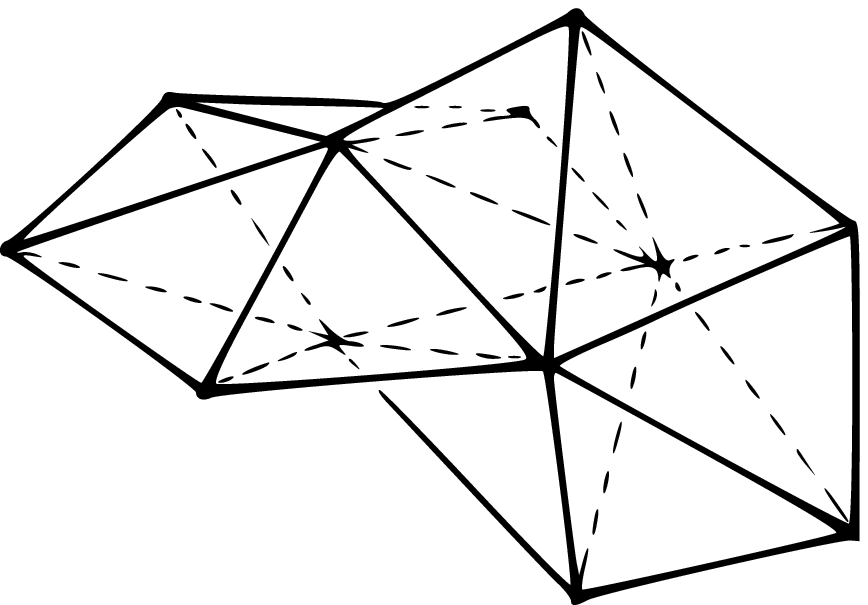, width=0.4\textwidth}
\end{figure}
\vskip -5mm \centerline{ \small Fig 1. A pentagonal Siamese
dipyramid}

\bigskip

It turns out that Siamese dipyramids have very interesting
deformability properties relied to the notion of flexibility.
Recall that a polyhedron with triangular faces is said to be {\it
flexible}, if it admits a continuous deformation, {\it a flexion},
which preserves both its combinatorial structure and its lengths
of edges but changes its spatial shape. The most known examples of
flexible polyhedra were constructed by Bricard, Connelly, Stefan,
see \cite{Connelly}, \cite[345–360]{Fuchs},
\cite[219–243]{Cromwell}.

\smallskip

As for the Siamese dipyramids, no one of them is flexible in the
sense of the classical theory of polyhedra. On the other hand,
their physical models made out of thin firm cardboard with hinged
joints of faces may behave similarly to physical models of
flexible polyhedra, i.e., these models may be unstable and admit
continuous deformations without any observable distortions
(ruptures, bucklings, breaks, extensions, contractions etc) of
faces but with significant variations in the spatial shapes. For
instance, this is the case for the equilateral pentagonal Siamese
dipyramids considered by M. Goldberg in \cite{Goldberg}. Thus the
flexibility of physical models of Siamese dipyramids contradicts
to the mathematical non-flexibility of these polyhedra.

\smallskip

The described phenomenon is called {\it the model flexibility}.
This phenomenon was observed earlier for various particular
polyhedra and it was claimed that the model flexibility can be
explained by slight non-destructive deformations of materials used
for producing models, cf \cite[p.224]{Cromwell}. However this
explanation was not confirmed completely by strong mathematical
reasonings.

\smallskip

Since the flexibility of polyhedra is based on the use of
flexions, i.e., continuous deformations which preserve both the
combinatorial structure and the metric structure of polyhedra,
then it is quite reasonable to assume that the model flexibility
of polyhedra may be explained with the help of {\it
 almost flexions}, i.e., continuous deformations which are {\it slight
modifications} of flexions. Namely, one can apply two kinds of
deformations of polyhedra: (1) continuous deformations which
preserve the combinatorial structure but variate slightly the
metric structure of polyhedra, (2) continuous deformations which
preserve the metric structure but change the combinatorial
structure. The both kinds of deformations were afore applied to
simulate the model flexibility of two particular polyhedra which
are the Alexandrov-Vladimirova star-like dipyramids
\cite{MilkaN}-\cite{Milka6} and the Jessen orthogonal icosahedron
\cite{GorkavyyKalinin}. We would strongly recommend for reading
the series of articles by A.D. Milka \cite{MilkaN}-\cite{Milka6},
where the model flexibility of Alexandrov-Vladimirova star-like
dipyramids was analyzed in the frames of the classical geometry
"in large", geometric theory of stability of shells developed by
A.V. Pogorelov, theory of dynamical systems, theory of
catastrophes etc.

\smallskip

In the present paper we discuss the model flexibility of Siamese
dipyramids by applying the approach based on the use of continuous
deformations preserving the combinatorial structure and varying
slightly the metric structure.

\smallskip

As the main result, we claim that for any $n\geq 3$ there exist $n$-gonal Siamese dipyramids which
admit continuous deformations preserving the combinatorial structure and such that negligible
relative variations in the length of edges produce significant relative variations in the spatial
shapes. For instance, for the equilateral pentagonal Siamese dipyramids there exists a continuous
deformation preserving the combinatorial structure and such that relative variations less than
0.004 in the lengths of edges generate relative variations greater than 5.9 in the heights
(distances between apexes) of the dipyramid. We stress that relative variations like 0.004 in the
lengths of edges may be treated as negligible and unobservable in physical models, that is why the
deformations mentioned above may be used to simulate the model flexibility of polyhedra in
question.

\smallskip

We hope that the proposed approach for explaining the model
flexibility can be adapted not only to Siamese dipyramids, but
also to any other model flexors.

\smallskip

The paper is organized as follows. In Sec. 1 we recall the definition of Siamese dipyramids,
describe their basic analytical and geometrical features and, in particular, demonstrate that no
one of the Siamese dipyramids is flexible. Particular examples are briefly presented in Sec. 2.
Next Sec. 3 concerns various approaches for representing analytically the whole family of Siamese
dipyramids. In Sec. 4 we discuss a rigidity map which represents a correspondence between extrinsic
(space shape) parameters and intrinsic (metric) parameters of Siamese dipyramids. Particular
attention is paid to topological properties of the rigidity map. Consequently we analyze the
existence of Siamese dipyramids with different space shapes but with the same metric properties
(the same lengths of edges). General continuous deformation of Siamese dipyramids are discussed
briefly in Sec.5. Next Sec.6 concerns deformability properties of equifacial Siamese dipyramids.
Particularly, we describe Siamese dipyramids which admit continuous deformations with relatively
small variations in the lengths of edges and hence may be used to illustrate the phenomenon of
model flexibility. Final Sec. 7 includes general conclusions and remarks.

\bigskip

\bigskip

\break

{\bf \large  1 Siamese dipyramids: definition and main features}

\medskip

Let us recall the definition of Siamese dipyramids. Fix a natural
$n\geq 3$, chose an arbitrary $l$ satisfying
\begin{eqnarray}\label{restriction}
\sin\frac{\pi}{2n} < \frac{l}{2} < \sin\frac{\pi}{n}
\end{eqnarray}
and consider an isosceles triangle $\Delta(l)$ with legs of length
$1$ and base of length $l$. Denote by  $\alpha$ the vertex angle
of $\Delta(l)$,
\begin{eqnarray}\label{alpha}
\sin\frac{\alpha}{2} = \frac{l}{2}.
\end{eqnarray}

Take $n$ copies of $\Delta(l)$ and stick them together to form a
triangulated ($n+2$)-gon $P(l)$ as shown in Fig. 2a. Due to
(\ref{restriction}), the polygon $P(l)$ is simple and non-convex,
it is bounded by $n$ edges of length $l$ and two edges of length
$1$. It is assumed that $P(l)$ may be folded along $n-1$ intrinsic
triangulating edges of length 1.

\begin{figure}[ht!]
    \centering
    \epsfig{file=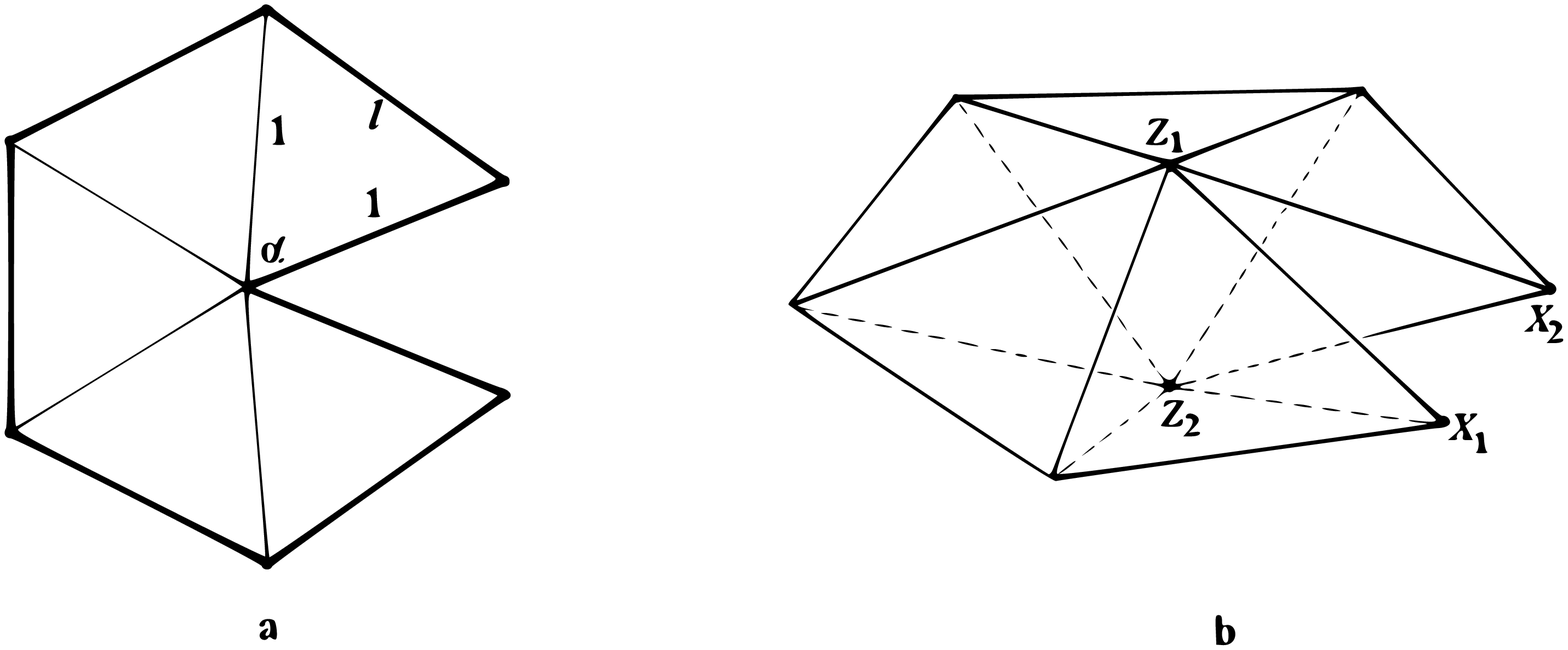, width=0.8\textwidth}
\end{figure}
\centerline{ \small Fig 2. a) A triangulated polygon
$P(l)$ with $n=5$, \,\, b) A Goldberg dipyramid $D(l)$ with $n=5$}

\bigskip

Next, take two copies of $P(l)$ and glue them together along the
edges of length $l$. This results in a dypiramid with boundary
which is called an {\it $n$-gonal Goldberg dipyramid}, see Fig.
2b. This dypiramid, $D(l)$, consists of $2n$ triangular faces
congruent to $\Delta(l)$, and its boundary $\partial D(l) =
Z_1X_1Z_2X_2$ consists of four segments of unit length.

The spatial shape of $D(l)$ is characterized by two quantities,
the {\it height} $\vert Z_1Z_2\vert = 2x$ and the {\it aperture}
$\vert X_1X_2\vert =2y$. Using elementary geometric arguments, one
can easily demonstrate that the following relation between $x$ and
$y$ holds true:
\begin{eqnarray}\label{relation}
y=\sqrt{1-x^2}\sin\left( n\,\arcsin\frac{l}{2\sqrt{1-x^2}}\right).
\end{eqnarray}
Each pair of positive $x$, $y$ satisfying (\ref{relation})
corresponds to a Godlberg dipyramid whose height and aperture are
equal to $x$ and $y$ respectively, we will denote it by $D(l;
x,y)$.

Actually, the relation (\ref{relation}) defines $y$ as a function
of $x$. If $x$ increases from $0$ to
$x_{max}=\sqrt{1-\frac{l^2}{4\sin^2\frac{\pi}{n}}}$ then $y$
decreases from $y_{max}=\sin \left( n\arcsin\frac{l}{2}\right)$ to
$0$. Hence for any fixed $n$ and $l$ satisfying (\ref{relation})
we obtain a one-parameter\footnote{Notice that each Goldberg
dipyramid in the family is well-defined by choosing either $x\in
[0, x_{\max}]$ or $y\in [0, y_{\max}]$, therefore $x$ or $y$ may
be used to parameterize the constructed family. Evidently, one can
apply any another parametric representation $x=x(t)$, $y=y(t)$,
$t\in [t_1,t_2]$, which satisfies (\ref{relation}) and guaranties
the positivity of $x$ and $y$.} continuous family of Goldberg
dipyramids which have the same combinatorial (symplicial)
structure  and the same lengths of edges but differ in the spatial
shapes. In other words, any Goldberg dipyramids is flexible in the
sense of the classical theory of polyhedra.

Finally, take two $n$-gonal Goldberg dipyramids, $D(l; x,y)$ and
$D(\tilde l; \tilde x,\tilde y)$, with possibly different $l$ and
$\tilde l$ fixed so that (\ref{restriction}) is satisfied, and
stick them together along their tetragonal boundaries so as to get
a closed polyhedron as shown in Fig. 1. This polyhedron is called
{\it an $n$-gonal Siamese dipyramid}, it consists of $2n$
isosceles triangular faces congruent to $\Delta(l)$ and $2n$
isosceles triangular faces congruent to $\Delta(\tilde l)$. If
$l=\tilde l$, then the faces are mutually congruent and the
Siamese dipyramid is referred to as {\it equifacial}.  If
$l=\tilde l=1$, then the faces are congruent to the equilateral
triangle with unit sides and the Siamese dipyramid is referred to
as {\it equilateral}.

Evidently, the described construction makes sense if and only if
the height and aperture of $D(l; x,y)$  are equal to the aperture
and height of $D(\tilde l; \tilde x,\tilde y)$ respectively, i.e.,
$\tilde y = x$, $\tilde x = y$. Expressing $\tilde y$ in terms of
$\tilde x$ and $\tilde l$ similarly to (\ref{relation}) and then
substituting $y = \tilde x$, $\tilde y = x$, we obtain the
following system:
\begin{eqnarray}\label{SiameseDef1}
\tilde x=\sqrt{1-x^2}\sin\left(
n\,\arcsin\frac{l}{2\sqrt{1-x^2}}\right),\\\label{SiameseDef2}
x=\sqrt{1-\tilde x^2}\sin\left( n\,\arcsin\frac{\tilde
l}{2\sqrt{1-\tilde x^2}}\right).
\end{eqnarray}
Thus each Siamese dipyramid, $S(l,\tilde l; x,\tilde x) = D(l;
x,\tilde x)\cup D(\tilde l; \tilde x,x)$, gives rise to a positive
solution $(x, \tilde x)$ of
(\ref{SiameseDef1})-(\ref{SiameseDef2}), and vice versa each
positive solution $(x, \tilde x)$ of
(\ref{SiameseDef1})-(\ref{SiameseDef2}) generates a well-defined
Siamese dipyramid. We will call $(x,\tilde x)$ the heights of
$S(l,\tilde l; x,\tilde x)$.

The solvability of (\ref{SiameseDef1})-(\ref{SiameseDef2}) with
respect to $x$ and $\tilde x$ depends on the values of $l$ and
$\tilde l$ viewed as auxiliary parameters for the system.
Actually, three different situations may happen:

1) the system (\ref{SiameseDef1})-(\ref{SiameseDef2}) has no
positive solutions, hence there no exist Siamese dipyramids with
faces congruent to $\Delta(l)$ and $\Delta(\tilde l)$;

2) the system (\ref{SiameseDef1})-(\ref{SiameseDef2}) has a unique
positive solution $(x,\tilde x)$, hence there exists a unique
Siamese dipyramid with faces congruent to $\Delta(l)$ and
$\Delta(\tilde l)$;

3) the system (\ref{SiameseDef1})-(\ref{SiameseDef2}) has at least
two different positive solutions, hence there exists at least two
Siamese dipyramids with faces congruent to $\Delta(l)$ and
$\Delta(\tilde l)$ -- these dipyramids have the same lengths of
edges but may differ in the spatial shapes.

Notice that the system (\ref{SiameseDef1})-(\ref{SiameseDef2})
possesses obvious symmetry properties. Namely, it remains
invariant under the interchange $\{x,l\}\leftrightarrow \{\tilde
x, \tilde l\}$, and this also affects its solvability.

If (\ref{SiameseDef1}) and (\ref{SiameseDef2}) were mutually
dependent for some choice of $l$ and $\tilde l$, then there would
exist a one-parameter continuous family of positive solutions of
(\ref{SiameseDef1})-(\ref{SiameseDef2}). This would result in a
one-parameter continuous family of non-congruent Siamese
dipyramids with the same lengths of edges. In this case one would
say that corresponding Siamese dipyramids are flexible similarly
to the well-known flexible polyhedra of Connelly and Steffan.

However (\ref{SiameseDef1}) and (\ref{SiameseDef2}) are not
mutually dependent, whatever $n$, $l$, $\tilde l$ are chosen. In
fact, these equations may be reduced to algebraical equations in
terms of $x$, $\tilde x$, $l$, $\tilde l$, whose concrete
expressions depend on $n$. If $n$ is odd, then (\ref{SiameseDef1})
is an algebraic equation of $n$-th order, which is of first order
with respect to $\tilde x$ and of $(n-1)$-th order with respect to
$x$. By symmetry, (\ref{SiameseDef2}) is an algebraic equation of
$n$-th order, which is of first order with respect to $x$ and of
$(n-1)$-th order with respect to $\tilde x$. If $n$ is even, then
(\ref{SiameseDef1}) is an algebraic equation of order $2n$, which
is of second order with respect to $\tilde x$ and of $2(n-1)$-th
order with respect to $x$. By symmetry, (\ref{SiameseDef2}) is an
algebraic equation of order $2n$, which is of second order with
respect to $x$ and of $2(n-1)$-th order with respect to $\tilde
x$. In both cases we may conclude that (\ref{SiameseDef1}) and
(\ref{SiameseDef2}) are mutually independent, therefore the system
in question may have a finite number of isolated solutions only,
whatever $n$, $l$, $\tilde l$ are chosen.

From the geometrical point of view, the following statements hold
true.

{\bf Claim 1.} {\it No one of the Siamese dipyramids is flexible.}

{\bf Claim 2.} {\it For any $n\geq 3$, $l$ and $\tilde l$
satisfying (\ref{restriction}), there exists at most finite number
of different  $n$-gonal Siamese dipyramids, whose faces are
congruent to $\Delta(l)$ and $\Delta(\tilde l)$.}

Let us emphasize that Siamese dipyramids mentioned in Claim 2 have
the same combinatorial structure and the same lengths of
corresponding edges. Such polyhedra are called {\it mutually
isomeric}, cf \cite[Ch.6]{Cromwell}.

\bigskip

\bigskip

{\bf \large 2 Concrete examples: original Siamese dipyramids}

\medskip

In order to illustrate the constructions described above, consider
some concrete examples including the case of equilateral
pentagonal Siamese dipyramids, which were discussed by Goldberg in
his original paper \cite{Goldberg}.

{\it Example 1.} Fix $n=5$ and let $l=1$, $\tilde l =1$. Then the
system (\ref{SiameseDef1})-(\ref{SiameseDef2}) rewrites as
follows:
\begin{eqnarray}\label{SiameseDef1Goldberg}
\tilde x = \frac{5x^4-5x^2+1}{2(1-x^2)^2},\\
\label{SiameseDef2Goldberg} x = \frac{5\tilde x^4-5\tilde
x^2+1}{2(1-\tilde x^2)^2}.
\end{eqnarray}
This system has exactly three different positive solutions: (a) $x\approx 0.07118$, $\tilde x
\approx 0.49237$, (b) $x\approx 0.32726$, $\tilde x \approx 0.32726$, (c) $x\approx 0.49237$,
$\tilde x \approx 0.07118$, see Fig. 3a. Hence, there exist three different equilateral pentagonal
Siamese dipyramids, see Fig. 4 below and Fig. 7 in \cite{Goldberg}. These three Siamese dipyramids
have the same combinatorial structure and the same lengths of edges, i.e., they are {\it isomeric},
but their spatial shapes are different.

\begin{figure}[ht!]
    \centering
    \epsfig{file=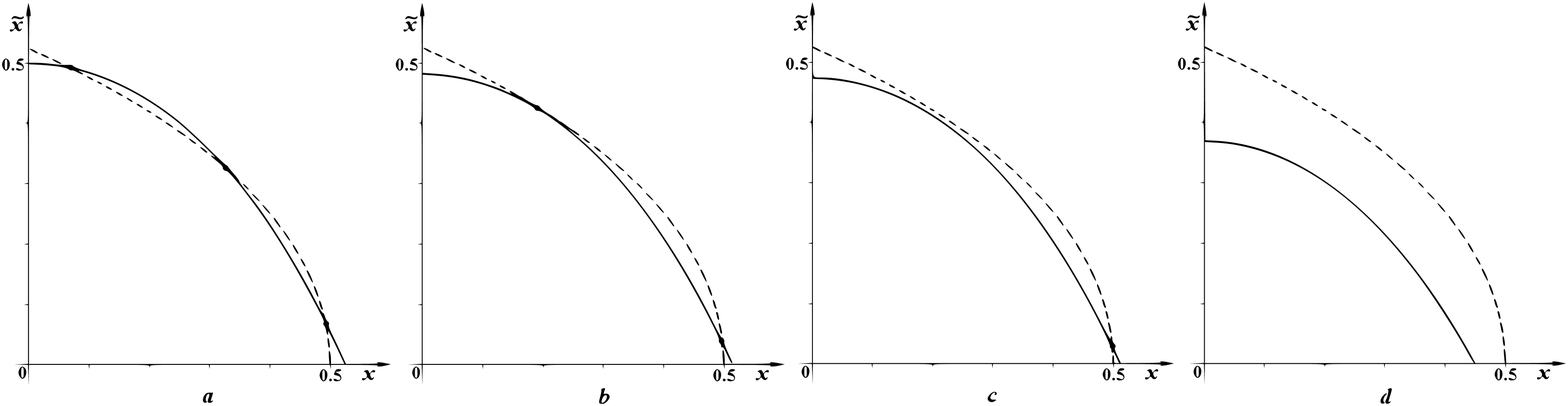, width=1\textwidth}
\end{figure}
\vskip -5mm \centerline{ \small Fig 3. Graphs of
(\ref{SiameseDef1})  and (\ref{SiameseDef2}) (dotted) for
different values of $l$ and $\tilde l$: }\centerline{ \small a)
$l=1$, $\tilde l =1$; b) $l\approx 1.0065$, $\tilde l =1$; c)
$l=1.01$, $\tilde l =1$; d) $l=1.05$, $\tilde l =1$}

\begin{figure}[ht!]
    \centering
    \epsfig{file=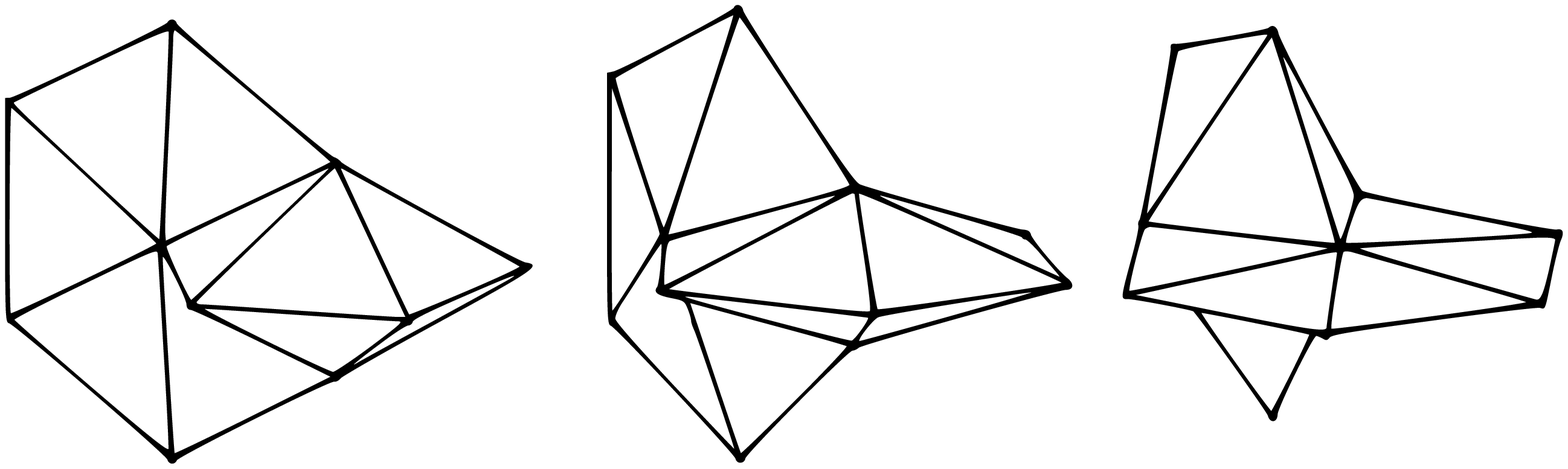, width=0.9\textwidth}
\end{figure}
\vskip -5mm \centerline{ \small Fig 4. Three different isomeric
pentagonal Siamese dipyramids with equilateral faces}

\bigskip

{\it Example 2.} Fix $n=5$ and let $l=1.01$, $\tilde l =1$. Then
the system (\ref{SiameseDef1})-(\ref{SiameseDef2})
has a unique positive solution  $x\approx 0.49888$, $\tilde x \approx 0.02721$, see Fig. 3c. Hence
for the given choice of $l$, $\tilde l$ there exists a unique Siamese dipyramid with faces
congruent to $\Delta(l)$ and $\Delta(\tilde l)$.

{\it Example 3.} Fix $n=5$ and let $l=1.05$, $\tilde l =1$. Then
the system (\ref{SiameseDef1})-(\ref{SiameseDef2})
has no positive solutions, see Fig. 3d. Hence for the given choice
of $l$, $\tilde l$ there no exist Siamese dipyramids with faces
congruent to $\Delta(l)$ and $\Delta(\tilde l)$.

{\it Example 4.} Fix $n=5$. Comparing Examples 1 and 2, one can recognize that there evidently
exists $l_0$ such that if one sets $l=l_0$, $\tilde l =1$, then the system
(\ref{SiameseDef1})-(\ref{SiameseDef2}) will has exactly two different positive solutions.
Actually, this is the case for $l_0\approx 1.0065$, and the solutions are (a) $x\approx 0.49773$,
$\tilde x \approx 0.03881$, (b) $x\approx 0.21214$, $\tilde x \approx 0.41334$, see Fig.3b. Hence
for the given choice of $l$, $\tilde l$ there exist two different isomeric Siamese dipyramids with
faces congruent to $\Delta(l)$ and $\Delta(\tilde l)$.

Notices, that Examples 1-3 may be viewed as generic, whereas
Example 4 is rather particular and may be qualified as a
transitional phase between Examples 1 and 2.


{\it Example 5.} An arbitrary $n\geq 3$ being fixed, consider an
arbitrary $l_0$ satisfying (\ref{restriction}) and set $l=\tilde l
= l_0$. Then the system (\ref{SiameseDef1})-(\ref{SiameseDef2})
has a positive solution such that $x=\tilde x$, this follows from
the symmetry of (\ref{SiameseDef1})-(\ref{SiameseDef2}) mentioned
above. Geometrically this means that for an arbitrary choice of
$n\geq 3$ and $l=\tilde l$ satisfying (\ref{restriction}) there
always exists an equifacial $n$-gonal Siamese dipyramid whose $4n$
isosceles triangular faces are congruent to the same triangle
$\Delta(l)=\Delta(\tilde l)$.

\bigskip

\bigskip

{\bf \large 3 Moduli space: analytical representation of Siamese
dipyramids}

\medskip

In Section 1  we derived the system
(\ref{SiameseDef1})-(\ref{SiameseDef2}) whose positive solutions
represent Siamese dipyramids. Namely, if one fixes arbitrary $l$
and $\tilde l$ satisfying (\ref{restriction}) and consider
(\ref{SiameseDef1})-(\ref{SiameseDef2}) as a system with respect
to $x$ and $\tilde x$, then for each positive solution $(x, \tilde
x)$ one can construct a well-defined $n$-gonal Syamese dipyramid
whose heights are equal to $x$ and $\tilde x$ and whose faces are
congruent to $\Delta (l)$ and $\Delta(\tilde l)$ respectively.
Therefore, it is natural to use
(\ref{SiameseDef1})-(\ref{SiameseDef2}) for  representing
analytically Siamese dipyramids. Let us discuss this idea in more
details.

Consider the first quadrant $\mathbb R^2_+ = \left\{ (x,\tilde
x)\vert x>0, \tilde x>0\right\}$ of the $(x, \tilde x)$-plane
$\mathbb R^2$. For every $l$ satisfying (\ref{restriction})
consider a curve $\gamma_l$ in $\mathbb R^2_+$ represented by
(\ref{SiameseDef1}). The family of curves $\gamma_l$, $l\in
(2\sin\frac{\pi}{2n},2\sin\frac{\pi}{n})$, foliates an open domain
$\Omega\in\mathbb R^2_+$ bounded by the curve
$\gamma_{2\sin\frac{\pi}{2n}}$ which is a limit curve for the
family, see Fig.5a.\footnote{Notice that another limit curve,
$\gamma_{2\sin\frac{\pi}{n}}$, degenerates to a point, the origin
$O$.} For each point $(x,\tilde x)\in \Omega$ there exists a well
defined $l\in (2\sin\frac{\pi}{2n},2\sin\frac{\pi}{n})$ such that
$\gamma_l$ goes through $(x,\tilde x)$. Actually, rewriting
(\ref{SiameseDef1}) in view of (\ref{restriction}), one gets the
following:
\begin{equation}\label{lll}
l = 2\sqrt{1-x^2}\sin\left(\frac{\pi}{n} -
\frac{1}{n}\arcsin\frac{\tilde x}{\sqrt{1-x^2}}\right).
\end{equation}
Geometrically this means that for each $(x,\tilde x)\in \Omega$
there exists a well defined $n$-gonal Goldberg dipyramid  whose
height and aperture are equal to $x$ and $\tilde x$ respectively.
Thus the whole family of $n$-gonal Goldberg dipyramids is well
represented by the points of the foliated domain
$\Omega\subset\mathbb R^2_+$.

Similarly, for every $\tilde l$ satisfying (\ref{restriction})
consider a curve $\widetilde \gamma_{\tilde l}$ in $\mathbb R^2_+$
represented by (\ref{SiameseDef2}). The family of curves
$\widetilde \gamma_{\tilde l}$, $\tilde l\in
(2\sin\frac{\pi}{2n},2\sin\frac{\pi}{n})$, foliates an open domain
$\widetilde \Omega\in\mathbb R^2_+$ bounded by the curve
$\widetilde \gamma_{2\sin\frac{\pi}{2n}}$; due to the symmetry of
(\ref{SiameseDef1}) and (\ref{SiameseDef2}), $\widetilde \Omega$
is symmetric to $\Omega$ with respect to the bisectrix of $\mathbb
R^2_+$. For each point $(x,\tilde x)\in \widetilde \Omega$ there
exists a well defined $\tilde l\in
(2\sin\frac{\pi}{2n},2\sin\frac{\pi}{n})$ such that $\widetilde
\gamma_{\tilde l}$ goes through $(x,\tilde x)$. Rewriting
(\ref{SiameseDef2}) in view of (\ref{restriction}), one gets the
following:
\begin{equation}\label{tilde l}
\tilde l = 2\sqrt{1-\tilde x^2}\sin\left(\frac{\pi}{n} -
\frac{1}{n}\arcsin\frac{x}{\sqrt{1- \tilde x^2}}\right).
\end{equation}
Geometrically this means that for each $(x,\tilde x)\in \widetilde
\Omega$ there exists a well defined $n$-gonal Goldberg dipyramid
whose aperture and height are equal to $x$ and $\tilde x$
respectively. Thus the whole family of Goldberg dipyramids is well
represented by the points of the foliated domain $\tilde
\Omega\subset\mathbb R^2_+$.

\begin{figure}[ht!]
    \centering
    \epsfig{file=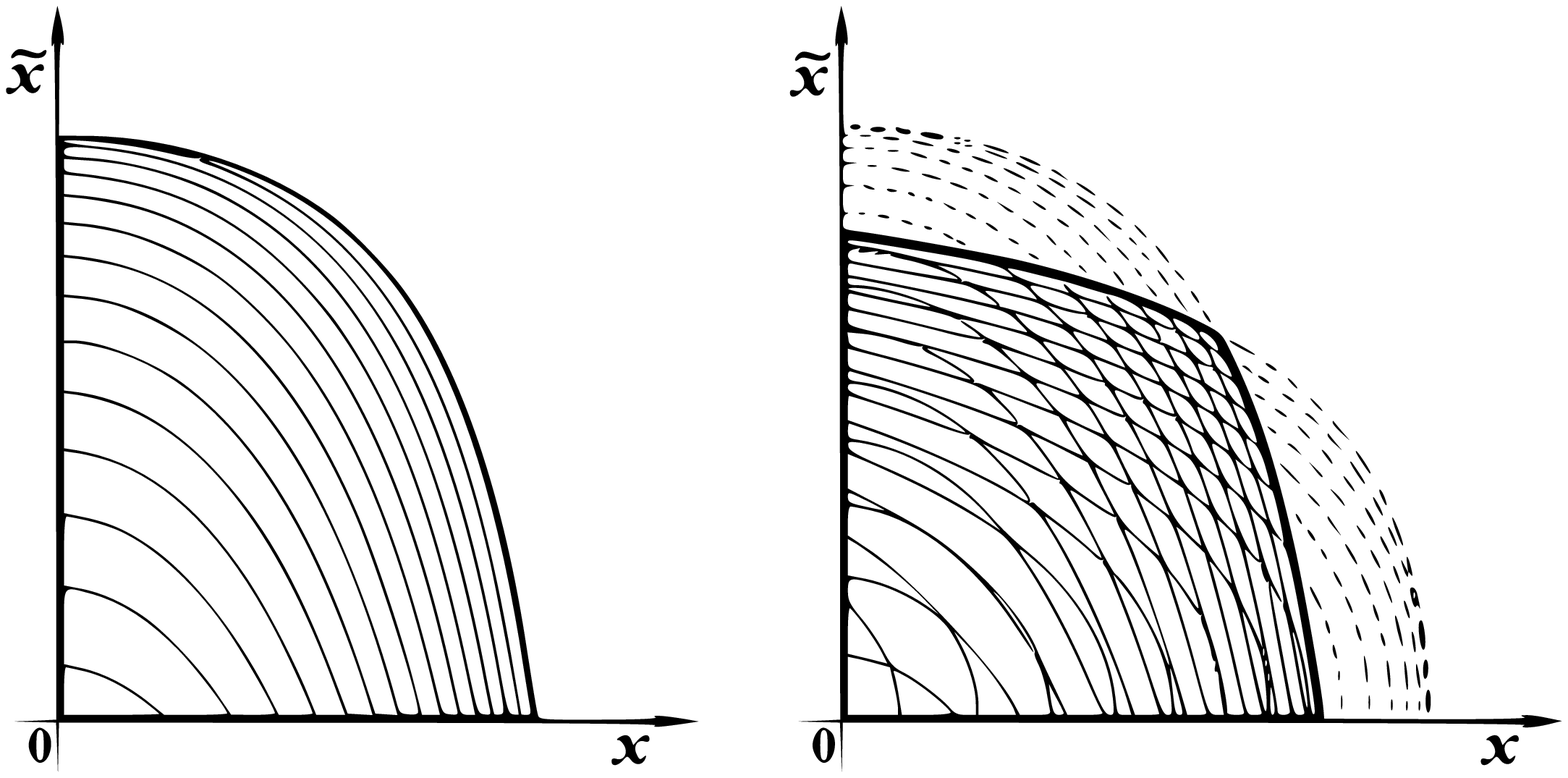, width=0.6\textwidth}
\end{figure}
\vskip -3mm \centerline{ \small Fig 5. The foliated domain
$\Omega$ (left) and the foliated domain $U$ obtained as the
intersection of foliated} \centerline{ \small   domains $\Omega$
and $\widetilde\Omega$ (right) }

\bigskip

Finally, consider the domain $U=\Omega \cap\widetilde\Omega$ in
$\mathbb R^2_+$, see Fig.5b. This open domain is foliated by
curves $\gamma_l\cap U$, $l\in
(2\sin\frac{\pi}{2n},2\sin\frac{\pi}{n})$, as well as by curves $
U\cap \widetilde\gamma_{\tilde l}$, $\tilde l\in
(2\sin\frac{\pi}{2n},2\sin\frac{\pi}{n})$. For each point
$(x,\tilde x)\in U$ there exist a well defined $n$-gonal Goldberg
dipyramid with height and aperture equal to $x$ and $\tilde x$
respectively and a well defined $n$-gonal Goldberg dipyramid with
aperture and height equal to $x$ and $\tilde x$ respectively.
Clearly, these Goldberg dipyramids may be glued together to form a
Siamese dipyramid with heights $x$ and $\tilde x$. Hence, the
following statement holds true.

{\bf Claim 3.} {\it For any $n\geq 3$ and $(x,\tilde x)\in U$
there exists a well-defined $n$-gonal Siamese dipyramid whose
heights are equal to $x$ and $\tilde x$. The lengths of edges of
this dipyramid are determined by (\ref{lll})-(\ref{tilde l}).}

Thus, the whole family of $n$-gonal Siamese dipyramids is well
represented by points of the domain $U\subset \mathbb R^2_+$.

Notice that this domain $U$ and its foliations in question depend
on $n$. However this dependence is not qualitative but rather
quantitative, since it reduces to a homothety after an appropriate
rescaling of the involved quantities $x$, $\tilde x$, $l$ and
$\tilde l$.

\bigskip

\bigskip

{\bf \large 4 Rigidity map: extrinsic vs intrinsic}

\medskip

Now consider a map $\varphi: U\subset \mathbb R^2_+ \rightarrow
\mathbb R^2$,
\begin{eqnarray*}
\varphi(x,\tilde x) = \left( l(x,\tilde x), \tilde l(x,\tilde x)
\right)
\end{eqnarray*}
represented by (\ref{lll})-(\ref{tilde l}), with $n\geq 3$ being
fixed. This is a well defined smooth map which assigns to the
heights of Siamese dipyramids the lengths of their edges. The
heights, $x$ and $\tilde x$, represent the spatial shapes of
Siamese dipyramids, whereas the lengths of edges, $l$ and $\tilde
l$, represent the metric structure of Siamese dipyramids. Thus,
$\varphi$ is a mapping between the extrinsic and intrinsic
geometries of $n$-gonal Siamese dipyramids, and therefore it may
be viewed as an analogue of the classical rigidity map arising in
the theory of polyhedra \cite{Gluck},
\cite{Connelly1}-\cite{Connelly2}.

The map $\varphi$ is defined in the domain $U\subset \mathbb
R^2_+$ of the $(x,\tilde x)$-plane. In view of
(\ref{restriction}), the image $\varphi(U)$ belongs to the square
$V$ in the $(l,\tilde l)$-plane,
\begin{eqnarray*}
V = \left\{ (l,\tilde l)\in \mathbb R^2\left\vert
\sin\frac{\pi}{2n} < \frac{l}{2} < \sin\frac{\pi}{n},
\sin\frac{\pi}{2n} < \frac{\tilde l}{2} < \sin\frac{\pi}{n}
\right. \right\}. \end{eqnarray*}

Examples in Section 3 demonstrate that the map $\varphi$ may be non-injective. Namely, if $n=5$,
then the pre-image of the point $(1,1)\in V$ consists of three points in $U$, which are
$(0.07118...,0.49237...)$, $(0.32726...,0.32726...)$ and $(0.49237..., 0.07118...)$; the pre-image
of the point $(1.01,1)\in V$ consists of a unique point $(0.49888...,0.02721...)\in U$; the
pre-image of the point $(1.05,1)\in V$ is empty; the pre-image of the point $(1.0065...,1)\in V$
consists of two points in $U$, which are $(0.49773...,0.03881...)$ and $(0.21214...,0.41334...)$.

In general, if two different points in $U$ have the same image in
$V$ under $\varphi$, then the $n$-gonal Siamese dipyramids
represented by these points are isomeric, i.e., they have the same
combinatorial structure and the same lengths of corresponding
edges. The converse is also true: if two $n$-gonal Siamese
dipyramids are isomeric, then they are represented by points in
$U$ which have the same image under $\varphi$. Thus, the
non-injectivity of the rigidity map $\varphi$ corresponds to the
isomericity of $n$-gonal Siamese dipyramids.

Hereinafter if $k$ points of $U$ are mapped by $\varphi$ to the
same point in $V$, then $k$ corresponding mutually isomeric
$n$-gonal Siamese dipyramids will be referred to as {\it
$k$-isomeric}. Clearly, any $1$-isomeric Siamese dipyramids is
well determined by assigning the lengths of its edges, but this is
not true for $k$-isomeric Siamese dipyramids with $k\geq 2$.

In order to understand the behavior of $\varphi: U\subset \mathbb
R^2_+ \rightarrow V\subset \mathbb R^2$, consider the boundary
$\partial U$ and the singular curve of $\varphi$ defined by
$\frac{\partial l}{\partial x} \frac{\partial \tilde l}{\partial
\tilde x} - \frac{\partial l}{\partial\tilde x}
\frac{\partial\tilde l}{\partial x} =0$, see Fig. 6a.

Their images under $\varphi$ generates a configuration of curves
in $V$ which bounds and partitions a domain in $V$ that is just
the image of $U$ under $\varphi$. It turns out that the
configuration of curves in question partitions $\varphi(U)\subset
V$ into five cells, see Fig. 6b. In its own turn, the full
pre-images of this configuration under $\varphi$ partitions $U$
into 9 cells, see Fig. 6c.

\begin{figure}[ht!]
    \centering
    \epsfig{file=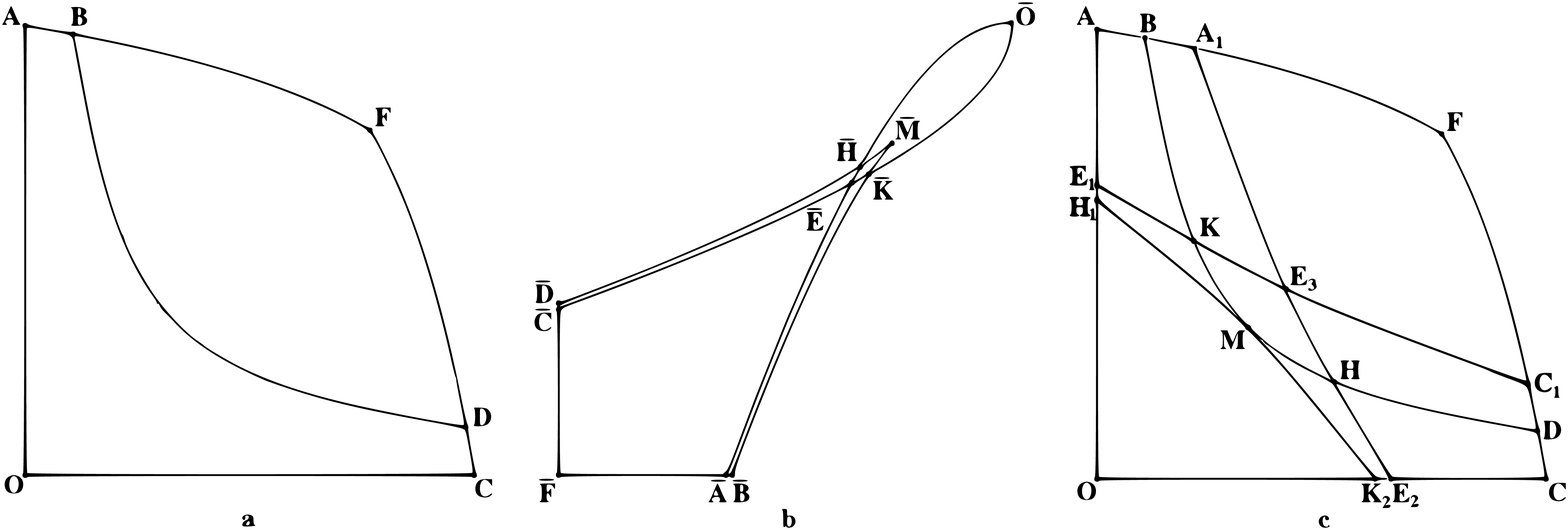, width=1\textwidth}
\end{figure}
\vskip -3mm \centerline{ \small Fig 6. a) Domain $U$, its boundary
$\partial U$ and the singular curve of $\varphi$} \centerline{
\small b) Image $\varphi(U)$ partitioned into 5 cells\,\, c)
Domain $U$ partitioned into 9 cells\footnote{ Characteristic
points in $U$ with the same image are labelled by the same
letters. The images of characteristic points are labelled by the
same letters as pre-images but overlined. } }
\bigskip

Consequently, the map $\varphi$ is cell-to-cell and its
restriction to any cell of $U$ is a diffeomorphism. Moreover, the
following holds true:

(i) There is one cell of $\varphi (U)$ which is covered by three
cells of $U$. In Fig. 6 this is the cell $\bar E\bar H\bar M\bar
K$ covered by three cells $E_1H_1MK_1$, $E_2H_2MK_2$,
$E_3H_2MK_1$. The points of this cells represents $3$-isomeric
Siamese dipyramids.

(ii) There are two cells of $\varphi (U)$ each of which is covered
by two cells of $U$. In Fig. 6 these are the cell $\bar C\bar
D\bar H\bar E$ covered by two cells $CDH_2E_2$, $C_1DH_2E_3$ and
the cell $\bar A\bar B\bar K\bar E$ covered by two cells
$ABK_1E_1$, $A_1BK_1E_3$. The points of this cells represents
$2$-isomeric Siamese dipyramids.

(iii) There are two cells of $\varphi (U)$ each of which is
covered by one cell of $U$. In Fig. 6 these are the cell $\bar
F\bar C\bar E\bar A$ covered by the cell $FC_1E_3A_1$ and the cell
$\bar O\bar H\bar M\bar K$ covered by the cell $OH_1MK_2$. The
points of this cells represents $1$-isomeric Siamese dipyramids.

It is important to emphasize that topologically the behaviour of
$\varphi$ is the same for any choice of $n\geq 3$, because by
changing the coordinates, $x^*=x$, $\tilde x^* = \tilde x$, $l^* =
\arcsin \frac{l}{2\sqrt{1-x^2}}$, $\tilde l^* = \arcsin
\frac{l}{2\sqrt{1-\tilde x^2}}$, one may easily conclude that
$\varphi$ rewritten in new coordinates depends on $n$
homotetically only. Thus the pictures drawn in Fig. 6 are
topologically the same for any $n\geq 3$, the only difference
concerns the concrete values of coordinates of characteristic
points in Fig. 6. These concrete values are pointed out in
Appendix 1.

From geometrical point of view we have the following.

{\bf Claim 4.} {\it For any $n\geq 3$, the set of all $n$-gonal
Siamese dipyramids contains:

-- an open subset of three-isomeric $n$-gonal Siamese dipyramids;

-- an open subset of two-isomeric $n$-gonal Siamese dipyramids;

-- an open subset of one-isomeric $n$-gonal Siamese dipyramids.

Besides, there no exist $k$-isomeric $n$-gonal Siamese dipyramids
with $k > 3$.}

Notice that equifacial Siamese dipyramids which correspond to points on the diagonal $\tilde l = l$
may be either three-isomeric (this is the case if $(l, l)$ is located in the segment $\bar E \bar
M$) or one-isomeric (this is the case if $l$ is sufficiently small so that $(l, l)$ belongs to the
segment $\bar F \bar E$ or sufficiently large so that $(l, l)$ belongs to the segment $\bar M \bar
O$).

\bigskip

\bigskip

{\bf \large 5 Continuous deformations of Siamese dipyramids: an
approach to the model flexibility}

\medskip

Now consider an arbitrary continuous curve $\Gamma: [0,T]\to U$.
Points of $U$ represent $n$-gonal Siamese dipyramids, therefore
$\Gamma$ generates a one-parameter continuous family of Siamese
dipyramids and may be interpreted as a continuous deformation of
Siamese dipyramids.

The variations of the heights of Siamese dipyramids under the
deformation are determined by the parametric representation of
$\Gamma$, i.e., $x=x(t)$, $\tilde x = \tilde x (t)$, this
describes the behavior of the spatial shapes of Siamese
dipyramids.

Substituting $x=x(t)$, $\tilde x = \tilde x (t)$ into
(\ref{lll})-(\ref{tilde l}), we obtain a parametric
representation, $l=l(t)$, $\tilde l = \tilde l (t)$, for the curve
$\varphi\circ\Gamma: [0,T]\to V$ which determines the variations
of the lengths of edges and hence controls the intrinsic geometry
of Siamese dipyramids under the deformation.

For instance, if $\Gamma$ was chosen so that $\varphi\circ\Gamma$
degenerated to a constant map, then $\Gamma$ would represent {\it
a flexion}, since the lengths of edges would remain constant under
the deformation. Clearly, this situation cann't happen since no
one of Siamese dipyramids is flexible, see Claim 1.

In this context, we will say that $\Gamma$ represents {\it an
almost flexion}, if $\varphi\circ\Gamma ([0,T])$ belongs to a
sufficiently small domain in $V$ so that $\Gamma$ represents a
deformation of Siamese dipyramids producing sufficiently small
relative variations in the lengths of edges. Analytically, a
subdomain in $V$ is viewed to be small if it belongs to a
rectangle centered at a point $(l_0,\tilde l_0)\in V$ and with
sides of lengths $2\varepsilon l_0$ and $2\varepsilon\tilde l_0$
parallel to the coordinate $l$- and $\tilde l$-axes, where
$\varepsilon >0$ is sufficiently small. If it is so then the
relative variations, $\delta l = \frac{\vert l-l_0\vert}{l_0}$ and
$\delta \tilde l = \frac{\vert \tilde l - \tilde l_0\vert}{\tilde
l_0}$, don't really exceed $\varepsilon$.\footnote{Implicitly we
make use here of the Chebyshev distance in $V$.}

Clearly, the definition above is meaningless until we fix a
specific value of $\varepsilon$. Particularly, if almost flexions
of polyhedra are aimed to simulate the flexibility of real world
polyhedral models, then the value of $\varepsilon$ is usually
chosen to be sufficiently small so that relative variations in the
lengths of edges may be qualified as
invisible/unobservable/inpalpable in real world models. For
example, in the frames of the geometric theory of stability of
shells developed by A.V. Pogorelov it is supposed that the value
of $\varepsilon$ should be comparable with $0.001$, i.e., $0.1\%$,
c.f. \cite[p.2]{PogGTS}.

\bigskip

\bigskip
\newpage

{\bf \large 6 Model flexibility of three-isomeric equifacial Siamese dipyramids}

\medskip

Let us focus on deformability properties of {\it three-isomeric
equifacial} Siamese dipyramids.

Namely, let $S_1$, $S_2$, $S_3$ be three different mutually
isomeric equifacial $n$-gonal Siamese dipyramids. They are
represented by three different points $P_1(x_1,\tilde x_1)$,
$P_2(x_2,\tilde x_2)$, $P_3(x_3,\tilde x_3$ in $U$ which a mapped
by $\varphi$ to the same point, a point $\bar P(l_0,\tilde l_0)\in
V$. This point is located in the cell $\bar H \bar M \bar K \bar
E$, whereas $P_1$, $P_2$ and $P_3$ belong to the cells
$E_1K_1MH_1$, $K_1E_3H_2M$, $H_2E_2K_2M$ respectively, recall the
cell decompositions shown in Fig. 6. The Siamese dipyramids are
assumed to be equifacial, hence we have $\tilde l_0 = l_0$ and
therefore $\tilde x_1 = x_3$, $\tilde x_2 = x_2$, $\tilde x_3 =
x_1$ by symmetry arguments.

Let $\Gamma : [0,T] \to U$ be a continuous curve which connects
successively $P_1$, $P_2$ and $P_3$. Then $\overline\Gamma =
\varphi\circ\Gamma : [0,T] \to V$ describes a double loop at $\bar
P \in V$. Since $P_2$ is separated from $P_1$ and $P_3$ by the
singular curve of $\varphi$, then the curve $\Gamma$ meets the
singular curve at least twice and hence the double loop
$\bar\Gamma$ has to meet at least twice the image of the singular
curve under $\varphi$.

The point $\bar P$ is said to be {\it admissible} if the straight
line $l + \tilde l = 2l_0$ through $\bar P$ in V meets the arcs
$\bar H \bar M$ and $\bar K \bar M$. Evidently, $\bar P$ is
admissible if it is sufficiently close to $\bar M$ so that the
inequality
\begin{eqnarray}\label{ineq}
\frac{l_{\bar H} + l_{\bar K}}{2} < l_0 < l_{\bar M}
\end{eqnarray} holds true, where $l_{\bar H}=\tilde l_{\bar H}$, $l_{\bar K}=\tilde l_{\bar K}$ and $l_{\bar
M}=\tilde l_{\bar M}$ stand for coordinates of $\bar H$, $\bar K$
and $\bar M$ respectively, see Fig. 7.

If $\bar P$ is admissible then we can specify $\Gamma$ so that the
points of this curve satisfy $l +\tilde l= 2l_0$ with $l$ and
$\tilde l$ expressed in terms of $x$ and $\tilde x$ by
(\ref{lll})-(\ref{tilde l}). In this case $\bar\Gamma$ describes a
twice covered segment $\hat H \hat K$ in $V$, where $\hat H$ and
$\hat K$ are the points where the straight line $l +\tilde l=
2l_0$ meets the arcs $\bar H \bar M$ and $\bar K \bar M$
respectively. When $\Gamma$ starts at $P_1$, goes through $P_2$
and ends at $P_3$, then $\bar\Gamma$ starts at $\bar P$, goes to
$\hat K$ , returns to $\bar P$, goes to $\hat H$ and finally
returns to $\bar P$, see Fig. 7.

\begin{figure}[ht!]
    \centering
    \epsfig{file=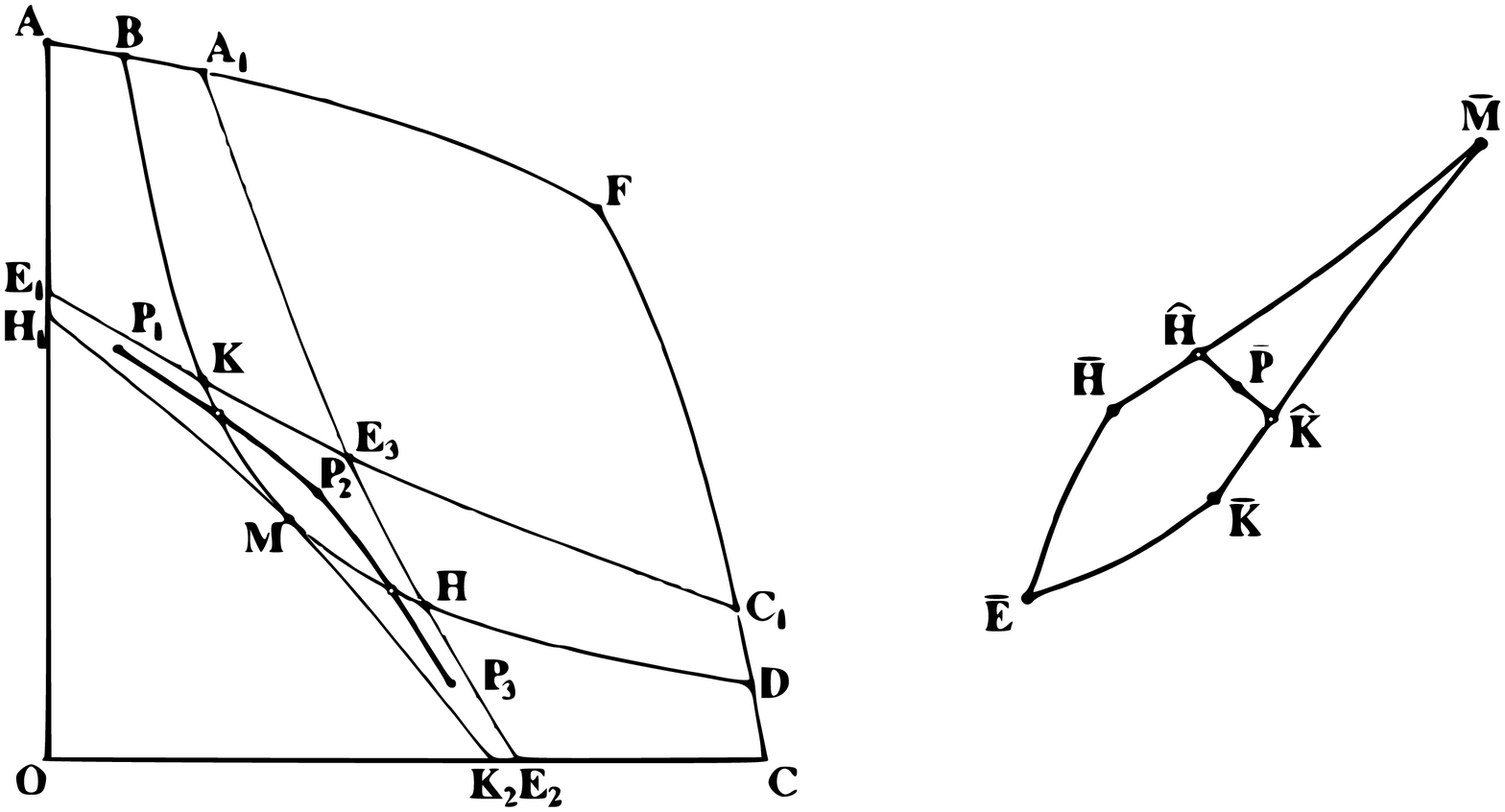, width=0.7\textwidth}
\end{figure}
\vskip -3mm \centerline{ \small Fig 7. Specified curve $\Gamma$ in
$U$ (left) and corresponding double loop $\overline\Gamma$ in
$\bar H\bar M\bar K \bar E \subset V$ (right)}

\bigskip

The continuous curve $\Gamma$ specified above represents a
continuous deformation of the Siamese dipyramids $S_1$, $S_2$ and
$S_3$ which will be referred to as {\it natural}. The relative
variations in the spatial shapes of Siamese dipyramids under this
deformation are estimated in terms of
\begin{eqnarray}\label{ddelta}
\delta_e = \max\left\{ \frac{\vert x_i - x_j\vert}{x_j}, 1\leq i, j\leq 3 \right\}.
\end{eqnarray}
The greater $\delta_e$ is, the more the spatial shapes of Siamese dipyramids in question differ one
from another.

As for the relative variations in the lengths of edges, $\delta l
= \frac{\vert l - l_0\vert}{l_0}$ and $\delta \tilde l =
\frac{\vert \tilde l - \tilde l_0\vert}{\tilde l_0}$, they may be
estimated via the length of the segment $\hat H\hat K$:
\begin{eqnarray}\label{intr}
\delta_i = \max\left\{ \delta l, \delta \tilde l \right\} =
\frac{\vert \hat H\hat K\vert}{2\sqrt{2}l_0}.
\end{eqnarray}
Notice that if $l_0$ increases and tends to $l_{\bar M}$ then
$\vert \hat H\hat K \vert$ decreases and tends to $0$, see Fig. 7.
Hence we get the following.

{\bf Claim 5.} {\it For any $n \geq 3$ and $\varepsilon > 0$,
there exists $l^* = l^*(\varepsilon, n)$ such that for any $l\in
(l^*, l_{\bar M} )$ there exist three different mutually isomeric
equifacial Siamese dipyramids with faces congruent to $\Delta(l)$
which may be connected by a continuous deformation in such a way
that the relative variations in the lengths of edges don’t exceed
$\varepsilon$.}

Consequently we may conclude that for any sufficiently small
$\varepsilon > 0$ there exist three different mutually isomeric
equifacial Siamese dipyramids which may be connected by a
continuous deformation with relative variations in the lengths of
edges less than $\varepsilon$. Evidently these Siamese dipyramids
provide us with an example of model flexibility: their physical
models would be unstable and behave like physical models of
flexible polyhedra.

For instance, set $l_0 = 2 \sin \frac{5\pi}{6n}$ whenever $n \geq
5$. This choice of $l_0$ means that $\alpha=\frac{5\pi}{3}$  and
hence $2\pi-n\alpha = \frac{\pi}{3}$, recall Fig. 1a.  Then it
turns out that $l_0$ satisfies (\ref{ineq}), c.f. Table 1, and
therefore $P(l_0,l_0)$ is admissible. Hence the corresponding
Siamese dipyramids $S_1$, $S_2$ and $S_3$ may be connected by a
natural continuous deformation. The value of $\delta_e$ turns out
to be sufficiently great, c.f. Table 1, hence the spatial shapes
of $S_1$, $S_2$ and $S_3$ are significantly different. On the
other hand, the value of $\delta_i$ turns out to be sufficiently
small, c.f. Table 1, and hence may be qualified as non-observable
in physical models of polyhedra.

$$
\begin{array}{cccccccccccc}
n  & 5 & 6 & 7 & 8 & 9 & 10 & 11 & 12&... \\
l_{\bar M}  & 1.02992  & 0.86634 & 0.74668 &
0.65565  & 0.58419 & 0.52667  & 0.47939  & 0.43986  &...\\
\frac{1}{2}\left(l_{\bar H}+l_{\bar K}\right) & 0.99438 & 0.83410
& 0.71767 & 0.62948 & 0.56045 & 0.50499 &  0.45947 & 0.42145 &...\\
l_0=2\sin\frac{5\pi}{6n} &  1  & 0.84523 & 0.73068  & 0.64287 &
 0.57360  & 0.51763  & 0.47151 & 0.43287  &... \\
\delta_i & 0.00394  & 0.00295 & 0.00240 & 0.00206 & 0.00184
& 0.00168 & 0.00157 & 0.00149   &...  \\
\delta_e & 5.91678  & 3.84229 & 3.03614 & 2.61522 & 2.36096
& 2.19322 & 2.07576 & 1.98986   &... \\
\end{array}
$$

\centerline{Table 1. Approximate values of characteristic parameters of Siamese dipyramids for
$5\leq n\leq 12$}

\bigskip

Thus we may conclude that the three-isomeric equifacial $n$-gonal
Siamese dipyramids with $l = \tilde l = 2 \sin \frac{5\pi}{6n}$,
$n \geq 5$, represent concrete examples for the phenomenon of
model flexibility. This conclusion may be easily verified by
producing and manipulating with cardboard models of the polyhedra
in question.

Notice that the cases of $n = 3$ or $4$ are rather particular, because $l_0 = 2 \sin
\frac{5\pi}{6n}$ does not satisfy the inequality (\ref{ineq}) that reads as $1.58292...< l_0 <
1.61379...$ if $n=3$ and $1.22698... < l_0 < 1.26433...$ if $n=4$. Hence another choice of $l_0$
has to be made to achieve the model flexibility in these cases. For instance, one may set $l_0
=1.6$ if $n=3$ and $l_0=1.25$ if $n=4$.

\bigskip

\bigskip

{\bf \large 7 Conclusion}

Evidently, model flexors show that the physical flexibility of practical polyhedral
con\-st\-ruc\-ti\-ons is not equivalent to the mathematical flexibility of corresponding polyhedra.
Siamese dipyramids as well as  Alexandrov-Vladimirova bipyramids and the Jessen icosahedron show
that the model flexibility of polyhedra is caused and explained by the presence of particular
continuous deformations such that negligibly small (almost invisible, unobservable) variations in
the lengths of edges provoke sufficiently large (visible, observable, palpable) variations in the
spatial shapes of polyhedra. This kind of deformations should be detected and analyzed for any
other polyhedral structures which pretend to illustrate  the phenomenon of model flexibility.

The discussed discrepancy between physical and mathematical notions of flexibility would be very
important for mechanics, architecture, technology, engineering, where mathematical simulations
using the theoretical flexibility of geometric objects (polyhedra, surfaces, etc) is applied to
ground the practical flexibility of real-world objects (shells, hinged plates, etc).

\bigskip

\bigskip


{\bf \large Appendix. Specific Siamese dipyramids and their characteristic features}

\bigskip

To complete the discussion above, let us describe particular Siamese dipyramids represented by
characteristic points shown in Fig. 6.

The point $A\in U$ with coordinates $(0,\tilde x_A)$ represents a
{\it degenerate} Siamese dipyramid $S_A$ with one height
vanishing, $x_A=0$. This means that one of two Goldberg pyramids
that compose $S_A$ collapses to a double covered planar polygon.
Moreover, $A\in\tilde \gamma_{2\sin\frac{\pi}{2n}}$ by definition,
hence for the lengths of edges of $S_A$ we have $\tilde l =
2\sin\frac{\pi}{2n}$. Consequently, the coordinates of $\bar A\in
V$ are $(l_A,2\sin\frac{\pi}{2n})$. The value of $\tilde x_A$ and
$l_A$ may be calculated numerically using
(\ref{SiameseDef1})-(\ref{SiameseDef2}) and
(\ref{lll})-(\ref{tilde l}) in view of $x_A=0$ and $\tilde
l_A=2\sin\frac{\pi}{2n}$, see Tables A1-A2.

The point $C$ is symmetric to $A$, i.e., $x_C=\tilde x_A$, $\tilde
x_C=0$. Moreover, the point $\bar C$ is symmetric to $\bar A$,
i.e., $l_C=\tilde l_A$, $\tilde l_C=l_A$. The corresponding
degenerate Siamese dipyramid $S_C$ is obtained from $S_A$ by
interchanging two Goldberg pyramids composing $S_A$.

The point $F\in U$ with coordinates $(x_F,\tilde x_F)$ belongs to
the curves $\gamma_{2\sin\frac{\pi}{2n}}$ and $\tilde
\gamma_{2\sin\frac{\pi}{2n}}$, hence it satisfy $\tilde x_F=x_F$
by symmetry. The corresponding point $\bar F\in V$ has coordinates
$(l_{\bar F}, \tilde l_{\bar F})=(2\sin\frac{\pi}{2n},
2\sin\frac{\pi}{2n})$. Therefore the Siamese dipyramid $S_F$
represented by $F$ is equifacial, all its faces are congruent to
the isosceles triangle $\Delta (2\sin\frac{\pi}{2n})$. Notice that
$S_F$ is the equifacial Siamese dipyramid with the smallest
lengths of edges $l=\tilde l$ as well as the equifacial Siamese
dipyramid with the greatest equal heights $x=\tilde x$, for any
fixed $n\geq 3$.

The point $O(0,0)\in U$ represents a degenerate Siamese dipyramid
$S_O$ with vanishing heights, $x_O=\tilde x_O =0$, therefore $S_O$
consists of two Goldberg pyramids collapsed to two double covered
$n$-gons. The lengths of corresponding edges of $S_O$ are equal to
$2\sin\frac{\pi}{n}$, this is the maximal possible value for $l$
and $\tilde l$. The corresponding point $\bar O\in V$ has
coordinates $l_{\bar O}=\tilde l_{\bar O} =2\sin\frac{\pi}{n}$.
Thus $S_O$ is the Siamese dipyramid with the maximal  possible
values of the lengths of edges, for any fixed $n\geq 3$.

The points of the singular curve of the rigidity map $\varphi$
represent shaky (infinitesimally flexible) Siamese dipyramids. The
point $B(x_B,\tilde x_B)$ is the point where the singular curve
meets the curve $\tilde \gamma_{2\sin\frac{\pi}{2n}}$, hence $B$
represents the shaky Siamese dipyramid, $S_B$, with $\tilde l =
2\sin\frac{\pi}{2n}$. The corresponding point $\bar B\in V$ has
coordinates $(l_{\bar B},\tilde l_{\bar B}) = (l_{\bar
B},2\sin\frac{\pi}{2n})$. Evidently, $S_B$ is the shaky Siamese
dipyramid with the maximal possible value of the height $\tilde x$
as well as the shaky Siamese dipyramid with the minimal possible
value of the length of edges $\tilde l$, for any fixed $n\geq 3$.
The value of $x_B$, $\tilde x_B$ and $l_{\bar B}$ may be
calculated numerically using
(\ref{SiameseDef1})-(\ref{SiameseDef2}), (\ref{lll})-(\ref{tilde
l}) and taking into account $\tilde l_{\bar B} =
2\sin\frac{\pi}{2n}$ and the singularity condition $\frac{\partial
l}{\partial x} \frac{\partial \tilde l}{\partial \tilde x} -
\frac{\partial l}{\partial\tilde x} \frac{\partial\tilde
l}{\partial x} =0$, see Tables A1-A2.

The point $D\in U$ is symmetric to $B$, i.e., $x_D=\tilde x_B$,
$\tilde x_D= x_B$. The corresponding point $\bar D\in V$ is
symmetric to $\bar B$, i.e., $l_{\bar D} = \tilde l_{\bar B}$,
$\tilde l_{\bar D} = l_{\bar B}$. Therefore the Siamese dypiramid
$S_D$ represented by $D$ may be obtained from $S_B$ by
interchanging two Goldberg pyramids composing $S_B$. Evidently,
$S_D$ is the shaky Siamese dipyramid with the maximal possible
value of the height $x$ as well as the shaky Siamese dipyramid
with the minimal possible value of the length of edges $l$, for
any fixed $n\geq 3$.

The point $M\in U$ is the symmetry point of the singular curve of
$\varphi$, its coordinates $x_M,tilde x_M)$ satisfy $\tilde
x_M=x_M$. The corresponding point $\bar M$ with coordinates
$(l_{\bar M},\tilde l_{\bar M})$ also satisfies $\tilde l_{\bar
M}= l_{\bar M}$. The Siamese dipyramid $S_M$ represented by $M$ is
the equifacial shaky Siamese dipyramid. Moreover, $S_M$ is the
shaky Siamese dipyramid with the maximal possible value of the
length of edges $l$ as well as the shaky Siamese dipyramid with
the maximal possible value of the length of edges $\tilde l$, for
any fixed $n\geq 3$. The value of $x_M$ and $l_{\bar M}$ may be
calculated numerically using
(\ref{SiameseDef1})-(\ref{SiameseDef2}), (\ref{lll})-(\ref{tilde
l}) and taking into account the symmetries $\tilde x_M=x_M$,
$\tilde l_{\bar M}= l_{\bar M}$ and  the singularity condition
$\frac{\partial l}{\partial x} \frac{\partial \tilde l}{\partial
\tilde x} - \frac{\partial l}{\partial\tilde x}
\frac{\partial\tilde l}{\partial x} =0$, see Tables A1-A2.

The rest of points, $A_1(x_{A_1},\tilde x_{A_1})$, $C_1(x_{C_1},\tilde x_{C_1})$, $H(x_H,\tilde
x_H)$, $K(x_K,\tilde x_K)$, $H_1(x_{H_1},\tilde x_{H_1})$, $K(x_{K_2},\tilde x_{K_2})$,
$E_1(x_{E_1},\tilde x_{E_1})$, $E(x_{E_2},\tilde x_{E_2})$, $E(x_{E_3},\tilde x_{E_3})$ in $U$ and
$\bar H (l_{\bar H},\tilde l_{\bar H})$, $\bar K(l_{\bar K},\tilde l_{\bar K})$, $\bar E(l_{\bar
E},\tilde l_{\bar E})$ in $V$ are interesting analytically rather than geometrically. They satisfy
by symmetry the following obvious relations: $x_{C_1}=\tilde x_{A_1}$, $\tilde x_{C_1}= x_{A_1}$,
$x_K = \tilde x_H$, $\tilde x_K=x_H$, $x_{K_2}=\tilde x_{H_1}$, $\tilde x_{K_2}=x_{H_1}$,
$x_{E_2}=\tilde x_{E_1}$, $\tilde x_{E_2}=x_{E_1}$, $\tilde x_{E_3} = x_{E_3}$ and $l_{\bar
K}=\tilde l_{\bar H}$, $\tilde l_{\bar K}=l_{\bar H}$, $\tilde l_{\bar E} = l_{\bar E}$, for
concrete values see Tables A1-A2.

$$
\begin{array}{cccccccccccc}
n & 3 & 4 & 5 & 6 & 7 & 8 & 9 & 10 & 11  & \infty\\
 x_A, \widetilde x_C & 0&0&0&0&0&0&0&0&0&0 \\
\widetilde x_A, x_C & .81649 & .84089 & .85065 & .85559 & .85847
& .86029 & .86152 & .86239 & .86303 & .86602 \\
x_F, \widetilde x_F & .64458 & .65063 & .65307 & .65432 & .65505 &
.65551 & .65582 & .65604 & .65621  &  .65697 \\
x_B, \widetilde x_D & .06239 & .08161 & .08979 & .09405 & .09655
& .09816 & .09924 & .10001 & .10058  & .10325 \\
\widetilde x_B, x_D & .80758 & .82710 & .83460 & .83833 & .84048 &
.84184 & .84275 & .84339 & .84387  &.84606 \\
x_M, \widetilde x_M & .21620 & .26650 & .28677 & .29710 & .30310 &
.30692 & .30949 & .31132 & .31266  & .31893 \\
x_H, \widetilde x_K & .34186  & .41946  & .45004  & .46550  &
.47446 & .48013  & .48396 & .48667  & .48865   & .49784 \\
\widetilde x_H, x_K & .13931  & .17044  & .18291  & .18921  &
.19286  & .19517 & .19673  & .19783 & .19864   &  .20238 \\
x_{E_1}, \widetilde x_{E_2} & 0&0&0&0&0&0&0&0&0&0 \\
\widetilde x_{E_1}, x_{E_2} & .42527 & .51922 & .55614 & .57469 &
.58540 & .59216 & .59672 & .59995 & .60231  & .61332  \\
x_{E_3}, \widetilde x_{E_3} & .27172 & .33219 & .35615 & .36824 &
.37524 & .37968 &
.38267 & .38478 & .38633  & .39357 \\
\end{array}
$$

\centerline{Table A1. Coordinates of characteristic points in $U$
for $3\leq n\leq 12$ and $n\to\infty$}

\bigskip

We would underline that the values of coordinates presented in
Table A1 depend on $n$ in such a way that this dependence is
essential for small values of $n$, but it almost disappears as
$n\to \infty$. This reveals a stabilization effect in the behavior
of the decomposed domain $U$ with respect to $n$ as $n\to \infty$.

$$
\begin{array}{cccccccccccc}
n & 3 & 4 & 5 & 6 & 7 & 8 & 9 & 10 & 11 \\
 l_{\bar A}, \tilde l_{\bar C}  & 1.33188 & 1.02082 &  .82441 &
.69045 & .59358 & .52037 & .46316 &
.41723 & .37956  \\
\tilde l_{\bar A}, l_{\bar C} &1. & .76536 &  .61803 & .51763 &
.44504 & .39018
& .34729 &  .31286  & .28462 \\
l_{\bar F}, \tilde l_{\bar F} & 1. & .76536  & .61803  & .51763 &
.44504 & .39018 &
.34729 & .31286  & .28462  \\
l_{\bar O}, \tilde l_{\bar O} & 1.73205   & 1.41421  & 1.17557  &
1 & .86776  & .76536  & .68404
& .61803  & .56346 0 \\
l_{\bar B}, \tilde l_{\bar D} & 1.33552 & 1.02599 &  .82965 &
.69535 & .59807  & .52447
& .47091 & .42067 & .38273  \\
\tilde l_{\bar B}, l_{\bar D}& 1. & .76536 & .61803 & .51763 &
.44504  &
.39018 & .34729 & .31286 & .28462 \\
l_{\bar M}, \tilde l_{\bar M} & 1.61379  & 1.26433  & 1.02992  &
.86634 & .74668 &
.65565  & .58419 & .52667  & .47939 \\
l_{\bar H}, \tilde l_{\bar K}&  1.57909  & 1.23259 & .98929 & .82944 & .71344 & .62565 & .55697 &
.50180 & .45654 \\
\tilde l_{\bar H}, l_{\bar K} & 1.58725 & 1.23954 & .99946 & .83876 & .72189& .63330 &
 .56393 & .50817 & .46240  \\
l_{\bar E} & 1.56761  & 1.20864  &  .97700  & .81837   &
.70353 & .61674  & .54890  & .49445  & .44979   \\
\end{array}
$$

\centerline{Table A2. Coordinates of characteristic points in $V$ for $3\leq n\leq 11$}

\bigskip

The values of coordinates presented in Table A2 depend on $n$ and
vanish as $n\to \infty$. However this dependence almost disappears
for great $n$, if one multiply the values of coordinates by $n$.
Hence if one replaces the decomposed domain $V$ by its
homotetically enlarged  copy $n V$, then one can easily observe a
stabilization effect in the behavior of $nV$ with respect to $n$
as $n\to \infty$.



\end{document}